\documentclass[a4paper]{amsart}

\usepackage{amssymb}
\usepackage{amsmath,amsthm}
\usepackage{enumerate,paralist}
\usepackage{amstext}
\usepackage{psfrag}
\usepackage{dsfont}
\usepackage[dvips]{epsfig}
\usepackage[english]{babel}
\usepackage{url}
\usepackage{graphicx}
\usepackage{dcolumn}
\usepackage{bm}
\usepackage{float}
\usepackage{hyperref}
\usepackage{color}
\usepackage{epstopdf}
\usepackage{cleveref}
\usepackage{cite}
\usepackage[svgnames]{xcolor}
\hypersetup{hidelinks,colorlinks=true,allcolors=DarkBlue}

\theoremstyle{plain}
\newtheorem{theorem}{Theorem}

\newtheorem{lemma}{Lemma}
\theoremstyle{definition}
\newtheorem{definition}{Definition}

\def\myref#1{(\ref{#1})}

\begin{document}
\frenchspacing

\title[Perturbation theory of observable linear systems]
{Perturbation theory \\ of observable linear systems}

\author{Aleksey Fedorov}
\address{
Institute for Problems in Mechanics, Russian Academy of Sciences \\ 119526, Vernadsky av., 101/1, Moscow, Russia \&
Russian Quantum Center \\ 143025 Novaya st. 100, Skolkovo, Moscow, Russia.}
 \email{akf@rqc.ru}

\author{Alexander Ovseevich}
\address{Institute for Problems in Mechanics, Russian Academy of Sciences \\
119526, Vernadsky av., 101/1, Moscow, Russia.} \email{ovseev@ipmnet.ru}

\maketitle

\begin{abstract}
The present work is motivated by the asymptotic control theory for a system of linear oscillators: the problem is to design
a common bounded scalar control for damping all oscillators in asymptotically minimal time. Motion of the system is
described in terms of a canonical system similar to that of the Pontryagin maximum principle. We consider evolution
equation for adjoint variables as a perturbed observable linear system. Due to the perturbation, the unobservable part of
the state trajectory cannot be recovered exactly. We estimate the recovering error  via the $L_1$-norm of perturbation.
This allows us to prove that the control makes the system  approach the equilibrium with a strictly positive speed.


\end{abstract}

\section{Introduction}
\subsection{Exact and approximate minimum time problem}
The subject of the present paper has grown out of study \cite{ovseev0, ovseev4}, of linear controllable dynamical system,
\begin{equation}\label{syst1}
    \dot{x}={A}x+{B}u,\quad x\in\mathbb{V}={\mathbb{R}}^{N}, \quad u\in\mathbb{U}={\mathbb{R}}, \quad|u|\leq1,
\end{equation}
where $A$ is a diagonalizable matrix with purely imaginary spectrum. The system is a standard model for control of
oscillations. In particular, if $N=2$ (one degree of freedom) we arrive at the classical problem of control of a single
oscillator by a bounded force, which is described in details in Ref. \cite{pont}.

Our basic  problem is the time-optimal damping of a given initial state of system (\ref{syst1}). Optimal trajectory is to
be found as the steepest descent in the direction of the gradient of the cost function, aka momentum.

The cost function and momenta are explicitly known in the case of one degree of freedom. For any greater number of
degrees of freedom the explicit description of the time-optimal damping is unknown, and, perhaps, does not exist.

\begin{definition}
The reachable set $\mathcal{D}(T)$ is the set of ends at time instant $T$ of all admissible trajectories of the system
starting at the given manifold at zero time.
\end{definition}

The level sets of the cost functions are boundaries of the reachable set of the system with respect to the backward time.
The direction of the gradient (momentum) is normal to the boundary of the reachable set.

Thus, the optimal control has the form
\begin{equation}\label{control_o}
    u(x)=-{\rm sign}\langle B,p(x) \rangle, \qquad p=\frac{\partial T}{\partial x}(x)
\end{equation}
where $p$ is  normal to the reachable set $\mathcal{D}(T(x))$ at point $x$ and angle brackets stand for the standard
scalar product in $\mathbb{R}^{N}$. We note that  $p=-\psi$, where $\psi$ is the standard adjoint variable of the
Pontryagin maximum principle. If the state $x$ is very far from the equilibrium, then the optimal damping time $T(x)$ is
large, and we need to know the behavior of $\mathcal{D}(T)$ as $T$ is large.

One of the basic results of \cite{ovseev} pertaining to system (\ref{syst1}) says that the reachable set $\mathcal{D}(T)$
equals asymptotically as $T\to\infty$ to the set $T\Omega,$ where $\Omega$ is a fixed convex body. E.g., if our system
describes a single oscillator, $\Omega$ is an ellipse in the plane $\mathbb R^2$. If the number of degrees of freedom is
greater, $\Omega$ does not possesses an elementary description, but still its support function has an explicit integral
representation.

We would like to study an easily implementable feedback control that asymptotically, as $x\to\infty$, works like the
optimal one. The basic idea of our control design is to substitute the set $T\Omega$ for $\mathcal{D}(T)$, and the
normals to this set for the momenta. If the phase vector $x$ lies in the boundary of $T\Omega$, then
\begin{equation}\label{approx3}
    x=T\frac{\partial {H}_\Omega}{\partial p}(p)
\end{equation} for a
momentum $p=p(x)$ and time $T$. Here ${H}_\Omega$ is the support function of $\Omega$ \cite{Hormander}. In the case at
hand the support function ${H}_{\Omega}$ is differentiable, and equation (\ref{approx3}) defines the direction of the
vector $p$ and the factor $T$ uniquely, because of the smoothness of the boundary of $\Omega$ \cite{ovseev0,ovseev2}.
Thus, the control is given by
\begin{equation}\label{approx55}
    u(x)=-{\rm sign}\langle{B,p(x)}\rangle.
\end{equation}

\subsection{Polar-like coordinate system}
Here we define a polar-like coordinate system, well suited for repre\-sen\-ta\-tion of motion under the control $u$.
Write the phase vector $x$ in the form $x=\rho\phi$, where $\rho>0$ and $\phi\in\omega=\partial\Omega$. In these
coordinates, equations of the motion have the form
\begin{equation}\label{T}
    \dot \rho=-\left|\left\langle\frac{\partial {\rho}}{\partial x},B\right\rangle\right|,\quad
    \dot\phi=A\phi+\frac{1}{\rho}\left(Bu+\phi \left|\left\langle\frac{\partial {\rho}}{\partial x},B\right\rangle\right|\right).
\end{equation}

It should be noted that $\rho=T$ and $\phi=\partial {H}_\Omega(p)/\partial p$ in terms of Eq. (\ref{approx3}). $\rho(x)$
is a norm of vector $x$, i.e., is a homogenous degree one convex function of $x$, which is strictly positive for
$x\neq0$.

The function $\rho$ is invariant under free (uncontrolled) motion of our system, i.e. $\langle{{\partial\rho}/{\partial
x},Ax}\rangle=0$. Moreover, we can show that the hessian ${\partial^{2}\rho}/{\partial x^{2}}$ is bounded on sphere
$|x|=1$.

An eikonal-type equation holds for the function $\rho=\rho(x)$
\begin{equation}\label{Euler}
    H_\Omega\left(p\right)=1,\quad p=\frac{\partial \rho}{\partial x}.
\end{equation}
It is ``dual'' to the equation $\rho({\partial {H_\Omega}}/{\partial p})=1$ of the surface $\omega$.

\subsection{Problem statement}
The following problem was addressed in Ref. \cite{ovseev0}: does control (\ref{approx55}) makes us approach zero with a
positive speed? It turns useful to measure distance to the target by means of the norm $\rho$. Precisely speaking,
suppose that $x(t)$ is a trajectory and $\rho_t=\rho(x(t))$, then, the question is: Is it true that
\begin{equation}\label{rate}
\rho_0-\rho_T\geq cT,
\end{equation}
where $c>0$? In other words, does our control make the norm $\rho$ approach zero with a strictly positive speed? By means of
theory presented below (Section \ref{Perturbation_theory}) we can give a positive answer provided that $\rho_t>M$ for
$t\in[0,T]$, where $M$ is a sufficiently large positive constant. We emphasize that all our results are related to the
motion very far from the origin.

In view of Eqs. (\ref{T})--(\ref{Euler}), we have to show that $L_1$-norm of $\left\langle p,B\right\rangle$ in the
interval $[0,T]$ is greater than $cT$. Time evolution of vector $p={\partial\rho}/{\partial x}$ is described as follows:
\begin{equation}\label{attractor_syst2p10}
    \dot p=-A^*p+\frac{\partial^2\rho}{\partial x^2}Bu=-A^*p+\widetilde Bu.
\end{equation}
Indeed, the total derivative is 
\begin{equation}
	\dot p=\frac{\partial^2\rho}{\partial x^2}(Ax+Bu). 
\end{equation}
Moreover, we have the identity,
\begin{equation}
    \left\langle{\frac{\partial\rho}{\partial x},Ax}\right\rangle=0,
\end{equation}
which expresses the invariance of the ``radius'' under free motion, where $u=0$. This latter identity implies that
\begin{equation}
    \frac{\partial^2\rho}{\partial x^2}Ax=-A^*p,
\end{equation}
and we arrive at Eq. (\ref{attractor_syst2p10}).

We note that the remainder $\widetilde Bu$ in Eq. (\ref{attractor_syst2p10}) is small provided that  the state vector $x$
moves at large distance from the origin, because ${\partial^{2}\rho}/{\partial x^{2}}$ is a homogeneous function of
degree $-1$ so that 
\begin{equation}
	\frac{\partial^2\rho}{\partial x^2}=O\left(\frac{1}{|x|}\right).
\end{equation}

Our idea is to interpret Eq. (\ref{attractor_syst2p10}) as the homogenous linear system
\begin{equation}\label{syst3}
    \dot p=-A^*p
\end{equation}
perturbed by $\widetilde Bu$, and consider $\left\langle p,B\right\rangle$ as a partial observation of the state vector
$p$. It allows us to embed our problem in the framework of the theory of observable linear systems \cite{kalman1,kalman2,kalman3}. The theory predicts that
the $L_1([0,T])$-norms of functions $p(t)$ and $\left\langle p(t),B\right\rangle$ are of the same order of magnitude, and
this leads to the required inequality: the $L_1([0,T])$-norm $\left\langle p,B\right\rangle$ is $\geq cT$.

\section{Perturbation theory}\label{Perturbation_theory}

The subject of the Kalman linear observation theory \cite{kalman1,kalman2,kalman3} is a linear time-invariant system,
$$
    \dot x=\mathcal{A}x, \quad y=\mathcal{C}x,
$$
which is observed, so that the vector $y$ is the observation result. Here $\mathcal{A}$ and $\mathcal{C}$ are constant
matrices.
\begin{definition}
    The system is said to be completely observable, if the knowledge of the curve $y(t)$ in an open time interval allows to recover $x(t)$ uniquely.
\end{definition}
We consider a perturbed situation where the observed vector has the same structure, but the vector $x$ satisfies the following
perturbed equation
\begin{equation}\label{perturbed}
    \dot x=\mathcal{A}x+f,\quad y=\mathcal{C}x.
\end{equation}
Then it is impossible to recover $x$ from knowledge of $y$ precisely, but, if the perturbation $f$ is small, we can do
this with a small error.

In quantitive terms, the error size is described by the following theorem:
\begin{theorem}\label{observation0}
Suppose that $\dot x=\mathcal{A}x,\,y=\mathcal{C}x$ is a completely observable time-invariant linear system. The
following a priori estimate holds for a solution $z$ of $\dot z=\mathcal{A}z+f$ in the interval $I$ of {\bf integer}
length
\begin{equation}\label{observ00}
    \int_I|z|dt\leq c\left(\int_I|\mathcal{C}z|dt+\int_I\left|f\right|dt\right),
\end{equation}
where the constant $c$ does not depend on the interval $I$.
\end{theorem}

{\bf Remark.} The condition that the interval $I$ has {\bf integer} length is not necessary. In order to make the
constant $c$ independent of $I$ it suffices to require that the length of $I$ is separated from zero. Without this
requirement, the corresponding statement is false.

\medskip

It is clear by summation over adjacent intervals of unit length that in order to prove Theorem \ref{observation0} it
suffices to consider the case $I=[0,1]$. We present two proofs of this theorem. The first one is easier, but less
constructive, the second one allows us in principle to find all involved constants efficiently.

\subsection{Direct proof}
The first proof of Theorem \ref{observation0} is based on the following Lemma:

\begin{lemma}\label{observation11}
Under the assumptions of Theorem \ref{observation0}, consider the map
\begin{equation}
     \Phi:z\mapsto [y,f]=[\mathcal{C}z,\dot z-\mathcal{A}z]
\end{equation}
from $\mathbb{W}=W^{1,1}\otimes\mathbb{R}^n$ to
$\mathbb{L}=\mathcal{L}_1\otimes\mathbb{R}^m\oplus\mathcal{L}_1\otimes\mathbb{R}^n$ and its image $L= \Phi(\mathbb{W})$.
Then the image $L$ of the map $\Phi$  is closed in $\mathbb{L}$.
\end{lemma}
Here $W^{n,1}$ is the Sobolev space of functions with $n$ integrable derivatives.

It is easy to derive Theorem \ref{observation0} from Lemma \ref{observation11}: The map $\Phi:\mathbb{W}\to L$ is a
continuous linear map. By Lemma \ref{observation11} the image $L$ is closed in $\mathbb{L}$. The observability condition
means that the kernel of the map $\Phi$ is zero. Hence one can apply the Banach inverse operator theorem and conclude
that
\begin{equation}\label{observ0}
    |z|_1\leq c(|\mathcal{C}z|_0+|\dot z-\mathcal{A}z|_0).
\end{equation}
Here $c$ is the norm of the inverse operator $\Phi^{-1}$, and
\begin{equation}
    |z|_n=\sum_{k=0}^n\int_0^1\left|\frac{\partial^k z}{\partial t^k}\right|dt
\end{equation}
is the standard Sobolev norm in $W^{n,1}([0,1])$. The conclusion of Theorem \ref{observation0} is an obvious relaxation
of inequality (\ref{observ0}).

To prove Lemma \ref{observation11}, we consider the subspace $M\subset L$ formed by vectors $\Phi(z)$ such that the
function $z$ vanishes at 0: $z(0)=0$. This is a closed subspace of $\mathbb{L}$, because the map $z\mapsto f=\dot
z-\mathcal{A}z$ defines an isomorphism $M\!\backsimeq\mathcal{L}_1\otimes\mathbb{R}^n$. Indeed, the Cauchy problem,
\begin{equation}
    \dot{z}=\mathcal{A}z+f,\quad z(0)=0,
\end{equation}
is correctly solvable. Another important subspace of $N\subset L$ is formed by  vectors $\Phi(z)$ such that $\dot
z-\mathcal{A}z=0$. It is also closed in  $\mathbb{L}$, because it is finite--dimensional $(\dim N=n)$. Since $L$ is a
direct sum of $M$ and $N$ it is  closed in $\mathbb{L}$. $\Box$

\subsection{Alternative proof of Theorem 1}
Another proof of Theorem \ref{observation0}, more constructive and more lengthy, is based on the following observation:
the change of parameters
\begin{equation}\label{brunovsky}
    \mathcal{A}\mapsto \mathcal{A}+\gamma \mathcal{C},\quad
    \mathcal{C}\mapsto \mathcal{C}, \quad
    \mathcal{A}\mapsto \delta{\mathcal{A}}\delta^{-1},\quad
    \mathcal{C}\mapsto{\mathcal{C}}\delta,
\end{equation}
where $\gamma$ is an arbitrary matrix, and $\delta$ is an arbitrary invertible matrix does not affect the validity of
Theorem, because it corresponds to substitutions $f\mapsto f+\gamma y$ and $z\mapsto \delta z$. In view of the Brunovsky
normal form \cite{brun} and the Kalman duality between controllability and observability, we may assume that the
observable system takes the form of a direct sum of systems with a scalar observation of the following form:
\begin{equation}\label{canonical}
\begin{array}{lcl}
    \dot z_1&=&z_2+f_1\\
    &\vdots&\\
    \dot z_{n-1}&=&z_n+f_{n-1}\\ \dot z_n&=&f_n
\end{array},
    \mbox{ and } y=z_1.
\end{equation}
These remarks reduce the proof of Theorem \ref{observation0} to the case of canonical system (\ref{canonical}) with the
scalar observation ($m=1$). We have to prove the following a priori estimate for solutions of system (\ref{canonical}):
\begin{equation}\label{observ000}
    |z|_0\leq c(|y|_0+|f|_0).
\end{equation}
In order to do so, we start with proving a relaxation
\begin{equation}\label{observ0001}
    |z|_0\leq c(|y|_0+|f|_n)
\end{equation}
of this inequality. It follows immediately from equations (\ref{canonical}) that
\begin{equation}\label{relation}
    \frac{\partial^n}{\partial t^n}y=\sum_{k=0}^n \frac{\partial^k}{\partial t^k}f_{n-k}.
\end{equation}
Denote the right-hand side (RHS) of identity (\ref{relation}) by $g$. It is clear that $|g|_0\leq C|f|_n$, and that $|z|_0\leq
C(|y|_n+|f|_n)$. Therefore, inequality (\ref{observ0001}) is implied by the following lemma:
\begin{lemma}\label{kolmogorov}Suppose that $y\in W^{n,1}$ and $\frac{\partial^n}{\partial t^n}y=g$. Then
\begin{equation}\label{observ0002}
    |y|_n\leq c(|y|_0+|g|_0).
\end{equation}
\end{lemma}
In other words, we have to estimate intermediate derivatives of $y$ by using estimates for 0-th derivative and $n$-th
derivative. This follows from the Kolmogorov-type inequality
\begin{equation}\label{Kolmogorov}
    \int_0^1\left|\frac{\partial y}{\partial t}\right|dt\ll\left(\int_0^1\left|\frac{\partial^ny}{\partial t^k}\right|dt\right)^{1/n}\left(\int_0^1\left|y\right|dt\right)^{(n-1)/n}
\end{equation}
proved in \cite{stein}.
Here $\ll$ is the Vinogradov symbol, meaning $O(\rm RHS)$.
Inequality (\ref{Kolmogorov}) implies:
\begin{equation}
    \int_0^1\left|\frac{\partial y}{\partial t}\right|dt\leq\epsilon\int_0^1\left|\frac{\partial^ny}{\partial t^n}\right|dt+\frac{C}{\epsilon}\left(\int_0^1\left|y\right|dt\right),
\end{equation}
where $C$ is a fixed constant, and $\epsilon>0$ is arbitrary.

In order to get from estimate (\ref{observ0001}) to (\ref{observ000}), we denote by $z_\phi,$ where $\phi\in L_1,$ the
solution of the Cauchy problem $\dot z_\phi=\mathcal{A} z_\phi+\phi$ with the initial condition $z_\phi(0)=z(0)$. It is
clear, e.g., from the Cauchy formula, that
\begin{equation}\label{1}
    |z-z_\phi|_0\ll|f-\phi|_0,
\end{equation}
and, therefore,
\begin{equation}\label{2}
    |\mathcal{C} z-\mathcal{C} z_\phi|_0\ll|f-\phi|_0.
\end{equation}
In view of (\ref{observ0001}), we obtain
\begin{equation}\label{4}
    |z_\phi|_0\leq c(|\mathcal{C} z_\phi|_0+|\phi|_n)
\end{equation}
Inequalities (\ref{1})--(\ref{4}) combined imply the a priori estimate
\begin{equation}\label{3}
|z|_0\ll |\mathcal{C} z|_0+\inf_\phi\{|f-\phi|_0+|\phi|_n\}
\end{equation}
for any $\phi\in L_1$.
Since $\displaystyle\inf_\phi \{|f-\phi|_0+|\phi|_n\}\leq|f|_0$  we arrive at estimate (\ref{observ000}). $\Box$

\medskip

{\bf Remark 1.} We note that  estimates (\ref{observ000}), (\ref{observ0001}), and (\ref{observ0002}) resemble basic
estimates in the $L_p$-theory of elliptic equations \cite{ADN}.

\medskip

{\bf Remark 2.} Note also that the canonical form (\ref{canonical}) allows us to give another proof of the crucial Lemma
\ref{observation11} in the first proof of our basic theorem. To do this, we invoke the following Lemma:
\begin{lemma}\label{observation111}
    There exists a pair of linear ordinary differential operators $P=P\left(\frac{\partial}{\partial t}\right)$ and $Q=Q\left(\frac{\partial}{\partial t}\right)$ with constant matrix coefficients such that $(y,f)\in L$ iff $Py=Qf$.
    Moreover, the degrees of polynomials $P,\,Q$ are $\leq n$.
\end{lemma}
This statement  immediately implies Lemma \ref{observation11}, because the condition $Py=Qf$ defines a closed subspace in
the space of (pairs of vector valued) distributions.  On the other hand, to prove Lemma \ref{observation111} it suffices
to consider only the case of scalar observable system (\ref{canonical}). But then, the relation $Py=Qf$ is given by
(\ref{relation}).

\section{Application}

We regard Eq. (\ref{attractor_syst2p10}) as a perturbed completely observable linear system (\ref{perturbed}), where the
phase vector $x=p$, matrices $\mathcal A=-A^*,\,\mathcal C=B^*$, observation $y=B^*p=\langle p,B\rangle$, and
perturbation $f={\widetilde B}u$. Assume that in the entire time interval $I$ of integer length $T$ the motion of the
state vector $x$ takes place within the domain $\rho(x)\geq C$. Then $|f|=O(1/C)$ in the entire interval. Moreover,
eikonal equation (\ref{Euler}) holds for $p$, and, therefore, $1\ll|p|$ and $T\ll\int_I|p|dt$.
The estimate of the main Theorem \ref{observation0} gives that
\begin{equation}
    \int_I|p|dt\ll \int_I|(p,B)|dt+\frac1C T,
\end{equation}
while the eikonal equation \myref{Euler} guarantees that $T\ll\int_I|p|dt$. By taking a sufficiently large constant
$C=C(A,B)$, we obtain that
\begin{equation}
    T\ll\int_I|(p,B)|dt.
\end{equation}
This is the inequality (\ref{rate}) in another notation. To be clear, we restate the result:

\begin{theorem}\label{observation5}
Suppose that  system \myref{syst1} moves from the level set $\rho=M$ to the level set $\rho=N$ under control
(\ref{approx55}), and $M,N\geq C(A,B)$, in the time interval of integer length $T$, where $C(A,B)$ is a (sufficiently
large) constant, depending only on parameters of system (\ref{syst1}). Then $T\leq c(M-N),$ where $c=c(A,B)$ is a
strictly positive constant.
\end{theorem}

\section*{Acknowledgements}

This work is supported by the Russian Foundation for Basic Research (grants 14-08-00606 and 14-01-00476). We thank the referee
and participants of seminar ``Problems in Differential Equations, Analysis and Control'' of the chair of General Control
Problems of the Lomonosov Moscow State University, especially A.V. Dmitruk, for helpful comments.

\end{document}